\documentclass[12pt]{article}
\usepackage[french,english]{babel}
\usepackage{amsmath}
\usepackage{amsfonts}
\usepackage[T1]{fontenc}
\usepackage[latin1]{inputenc}
\usepackage{graphicx}
\def\as={\,\stackrel{a.s.}{=}\,}
\newdimen\AAdi
\newbox\AAbo

\def\AAk#1#2{\setbox\AAbo=\hbox{#2}\AAdi=\wd\AAbo\kern#1\AAdi{}}
\def\AAr#1#2#3{\setbox\AAbo=\hbox{#2}\AAdi=\ht\AAbo\raise#1\AAdi\hbox{#3}}

\newcommand{\R}{{\mathbb{R}}}

\newtheorem{lmm}{Lemma}

\newtheorem{thrm}{Theorem}
\newtheorem{coro}{Corollary}

\begin{document}

\title{ Asymptotic normality of  recursive  estimators under strong mixing conditions }

\author{
{\large  Amiri Aboubacar$^{}$
\footnote{
Universit\'e Lille Nord de France, Universit\'e Lille 3,   Laboratoire EQUIPPE
EA 4018, Villeneuve d'Ascq, France. aboubacar.amiri@univ-lille3.fr.
} 
} \\
}

\maketitle
\thispagestyle{empty}

\vspace{1.1cm}

\renewcommand{\baselinestretch}{0.9}

\def\baselinestretch{1.3}
\pagestyle{plain}
\setcounter{page}{1}
\normalsize
\setcounter{footnote}{0}

\begin{abstract}
The  main purpose of this paper is to  estimate the regression function  by using a recursive nonparametric kernel approach.  We derive the  asymptotic normality for a general class  of recursive kernel estimate of the regression function, under  strong mixing conditions.  Our purpose is  to extend the work of Roussas and Tran \cite{Roussas-Tran}  concerning the Devroye-Wagner  estimate.

\end{abstract}


\section{Introduction}\label{intro}
In this paper we consider  nonparametric sequential estimation of a regression functional, for dependent observations. Regression function estimation is an important problem in data analysis and remains a subject of hight interest, which covers many applied fields  such  as prediction, econometrics,  decision theory,  classification, communications and control systems.   The literature on this topic is still growing and some relevant work on this subject include the monographs by  Prakasa Rao\cite{prakassa-rao1}, Gy\"orfi et al. \cite{GyorfiHardle} and Yoshihara \cite{Yoshihara}, while more recent results  are presented  in, for example, the books by Gy\"orfi et al. \cite{GyorfiKholer} and Bosq and Blanke \cite{bosq-blanke}.
 Sequential estimation is achieved with the use of recursive estimators, typically kernel ones,  and the purpose of this paper is to study a certain class of them.  An estimator is said to be  `recursive' if its value calculated from the first $n$ observations, say $f_n$, is only a function of $f_{n-1}$  and the  $n^{th}$ observation.  In this way, the estimator can be updated with each new observation added to the database.  This  recursive property  is clearly useful in sequential investigations  and also for a fairly large sample size, since addition of a new observation means the non-recursive estimators must be entirely recomputed. Besides, we are required to store extensive data in order to re-calculate them. 
 
  The first kernel recursive regression estimator was introduced by  Ahmad and Lin \cite{Ahmad-lin} taking the form 
$$ r_n^{AL}({x}) :=\frac{\sum\limits_{i=1}^nY_iK\left(\frac{x-X_i}{h_i}\right)}
{\sum\limits_{i=1}^nK\left(\frac{x-X_i}{h_i}\right)},$$ 
 which  is a  recursive version of  the Nadaraya-Watson   estimate.  Also  Devroye and Wagner \cite{dev-wag} propose the recursive estimator of the form 
 $$r_n^{DW}({x}) :=\frac{\sum\limits_{i=1}^n\frac{Y_i}{h_i}K\left(\frac{x-X_i}{h_i}\right)}
{\sum\limits_{i=1}^n\frac{1}{h_i}K\left(\frac{x-X_i}{h_i}\right)}. $$
In the literature  $r_n^{AL}({x})$ and   $r_n^{DW}({x})$  are  respectively  the so-called recursive and semi-recursive estimators.  Various results on the latter estimators  were established in an independent and identically distributed (i.i.d.) case, by many authors, we cite, among many others, Ahmad and Lin \cite{Ahmad-lin}, Devroye \cite{devroye81}, Greblicki and Pawlak \cite{Grebilicki-pawlak}, Krzyzak \cite{KRZ} and Walk \cite{walk}.   In the dependent case, the majority of works  are focused to Devroye-Wagner estimate. 
In a context of strong mixing processes, Roussas \cite{Roussas90} gave the uniform almost sure convergence for $r_n^{DW}(x)$, and Roussas and Tran \cite{Roussas-Tran} showed its asymptotic normality. Under $\varphi$-mixing conditions, Qin \cite{Qin} showed the asymptotic normality of $r_n^{DW}(x)$, and Wang and Liang \cite{Wang} studied the almost uniform convergence for truncated versions of $r_n^{DW}(x)$ and $r_n^{AL}(x)$ in the same context.
It should be noted that, unlike to the iid case,  more results are only  obtained for $r_n^{DW}(x)$ in  dependent  case.  In particular no asymptotic normality has so far been established for $r_n^{AL}(x)$   in this context. Also we remark that, the approach utilized by Roussas and Tran  \cite{Roussas-Tran} to establishing the asymptotic normality of  $r_n^{DW}(x)$ cannot be generalized step by step to $r_n^{AL}(x)$.  Indeed,  the adaptation  of  their proof to $r_n^{AL}(x)$,  needs  to suppose that the sequence  $\frac{1}{n}\sum_{i=1}^n(h_i/h_n)^{2d}$ converges to a  finite limit, for the study  of a few   covariance terms.  The earlier condition is not satisfied by the popular choice 
$h_n=cn^{-\frac{1}{d+4}}$ for $d>3$.  Also their proof uses the fact that  for all $i=1,\dots,n$ $h_n<h_i$, while the same approach applied to  $r_n^{AL}(x)$,  leads to assume that $h_n>h_i$, which contradicts the optimal  choice of $h_n$.  \\
This paper deals with an extension  of the work of Roussas and Tran \cite{Roussas-Tran}  to  the general family  of recursive estimators  introduced by  Amiri \cite{Amiri},  whose $r_n^{DW}(x)$ and $r_n^{AL}(x)$  are special cases.   
The paper is organized as follows.  In the next section, we present our main assumptions and the result for regression estimation. The proof of the  main result is postponed until  Section \ref{preuves}.
  \section{ Sequential regression estimation}
\subsection{Notation and assumptions}
Let $\left\{(X_t,Y_t), t\in\mathbb{N}\right\}$
be  a sequence of random variables on probability space $\left(\Omega,{\cal F}, P\right)$, taking values in $\R^d\times\R^{d'}$ ($d\geq1,d'\geq1$), and having probability density function $f_{(X,Y)}$
with respect to the  Lebesgue mesure. Assume that  $m$ is a  Borelian  function  on $\R^{d'}$ into  $\R$ such that $\omega\mapsto m^2\left( Y_t(\omega)\right) $ is  $P$-integrable, and define  the regression function as
 $$r(x):=\left\{\begin{array}{l}
\text{E}\left(m\left( Y_0\right) |X_0=x \right)=\dfrac{\int_{\R^{d'}}m(y)f_{(X,Y)}(x,y)dy}{f(x)}:=\frac{\varphi(x)}{f(x)},\text{ if } f(x)>0\\
\text{E}m(Y_0), \text{ if } f(x)=0, \end{array}\right.$$
where $f$ is the probability density function of $X_0$. Note that the transformation $m$ is chosen by the statistician leading to multiple choices of estimation. Typical examples  of $m$ are  identity and polynomial functions to estimate respectively the usual regression and the conditional moments. \\

 Throughout the paper we suppose that $f, \varphi\in C_d^2(b),$ where $C_d^2(b)$ denotes the set of twice-differentiable functions,  with bounded second derivative.  This condition is  classical in the area of nonparametric  estimation and  has been  used, among others,  by Roussas and Tran \cite{Roussas-Tran},  Bosq and Blanke  \cite{bosq-blanke}.   

To estimate the functional $r(x),$  we consider the general family  of kernel  regression estimators introduced in Amiri \cite{Amiri},    defined by
\begin{equation}\label{reg}
r_n^\ell(x):=\frac{\sum\limits_{i=1}^n\dfrac{m\left(Y_i\right)}{h_i^{d\ell}}K\left( \dfrac{x-X_i}{h_i}\right)}{\sum\limits_{i=1}^n\dfrac{1}{h_i^{d\ell}}K\left( \dfrac{x-X_i}{h_i}\right)}, \end{equation}
which can be computed recursively by
$$r_{n}^\ell(x)=\frac{\left(\sum\limits_{i=1}^{n-1}h_i^{d(1-\ell)}\right)\varphi_{n-1}^\ell(x)+\left(\sum\limits_{i=1}^{n}h_i^{d(1-\ell)}\right)m(Y_{n})K_{n}^\ell\left(x-X_{n} \right)}
{\left(\sum\limits_{i=1}^{n-1}h_i^{d(1-\ell)}\right)f_{n-1}^\ell(x)+\left(\sum\limits_{i=1}^{n}h_i^{d(1-\ell)}\right)K_{n}^\ell\left(x-X_{n} \right)},$$
where $$\varphi_n^\ell(x):=\dfrac{1}{\sum\limits_{i=1}^nh_i^{d(1-\ell)}}\sum\limits_{i=1}^n\dfrac{m(Y_i)}{h_i^{d\ell}}K\left( \dfrac{x-X_i}{h_i}\right),~f_n^\ell(x):=\dfrac{1}{\sum\limits_{i=1}^nh_i^{d(1-\ell)}}\sum\limits_{i=1}^n\dfrac{1}{h_i^{d\ell}}K\left( \dfrac{x-X_i}{h_i}\right),$$ and $K_i^{\ell}(\cdot):=\frac{1}{h_i^{d\ell}\sum\limits_{j=1}^ih_j^{d(1-\ell)}}K\left(\frac{\cdot}{h_i}\right)$.  
Our class of estimates   includes the  popular kernel recursive estimators 
$r_n^{AL}({x})$ and   $r_n^{DW}({x})$,  corresponding to the cases $\ell=0$ and $\ell=1$, respectively.\\
 
At this point,  we make some assumptions  and give the main theorem. Throughout this paper the kernel $K$ is assumed to satisfy the following conditions.\\

\textbf{\underline{Assumption $\textbf{H1}$.}}
\begin{description}
\item[(i)] $K:\mathbb{R}^d\mapsto\mathbb{R}$ is bounded, symmetric and positive function such that  $\int_{\R^d}K(t)dt=1;$
\item[(ii)] $\lim\limits_{\left\|x\right\|\rightarrow{+\infty}}\left\|x\right\|^dK(x)=0;$
\item[(iii)] $\int_{\mathbb{R}^d}|v_iv_j|K(v)dv<\infty,~~i,~j=1,\ldots,d.$\\
\end{description}

Assume the sequence   $h_n$ satisfies   the  following conditions. \\

{\textbf{\underline{Assumption  $\textbf{H2}$.}}}
\begin{description}
  \item[(i)] $h_n\downarrow0, nh_n^{d+2}\rightarrow\infty; $
  \item[(ii)] For all $r\in]-\infty,~~d+2]$,   $B_{n,r}:=\frac{1}{n}\sum\limits_{i=1}^n\left(\frac{h_i}{h_n} \right)^{r} \rightarrow\beta_{r}>0\text{ as }n\rightarrow\infty;$
  \item[(iii)]  For each  sequence of integers $u_n$ and $v_n$ such that  $u_n\sim v_n$, then  $h_{u_n}\sim h_{v_n}.$ \footnote{If  $a_n$ and $b_n$ are two real  sequences,   $a_n\sim b_n$ means that the ratio $a_n/b_n$ converges 1} \\
 \end{description}

{\textbf{\underline{Assumption  $\textbf{H3}$.}}}
\begin{description}
\item[(i)] The  process  $(X_t)$ is  $\alpha$-mixing with  $$\alpha_X(k)\leq\gamma k^{-\rho},~~k\geq 1, \gamma>0\text{ and }\rho>\max\left(2,\frac{d+2}{2}\right);$$
\item[(ii)] For each couple $(s,t),~~s\neq t,$ the random vector  $(X_s,X_t)$
admits  a probability density function  $f_{(X_s,X_t)}$  such that
 $\sup\limits_{\mid s-t\mid\geq1}\parallel g_{s,t}\parallel_{\infty}<\infty,\text{ where } g_{s,t}(\cdot,\cdot):=f_{(X_s,X_t)}(\cdot,\cdot)-f(\cdot) f(\cdot).$
\end{description}

\underline{\textbf{Assumption  $\textbf{H4}$}.}
\begin{description}
\item[(i)] The function $\text{E}\left(m^2(Y)|X_0=\cdot\right)f(\cdot)$ is both continuous and bounded away from zero at $x$;
\item[(ii)]$\text{There exist }\lambda>0,\theta>0 \text{ such that }\text{E}\exp(\lambda |m(Y_0)|^\theta)<\infty;$
\item[(iii)] For each   $k\neq k',$ the random vector $\left(X_k,Y_k,X_{k'},Y_{k'}\right)$  admits a probability density function   $f_{\left(X_k,Y_k,X_{k'},Y_{k'}\right)}$, such that $\sup\limits_{|k-k'|\geq1}\sup\limits_{(s,t)\in\R^{2d}}\int_{\mathbb{R}^{d'}}\int_{\mathbb{R}^{d'}}\left| G_{k,k'}\left(s,u,t,v\right)\right|dudv<\infty,$ where\\ $G_{k,k'}\left(\cdot,\cdot,\cdot,\cdot\right)=f_{\left(X_k,Y_k,X_{k'},Y_{k'}\right)}\left(\cdot,\cdot,\cdot,\cdot\right)-f_{(X,Y)}\left(\cdot,\cdot\right)f_{(X,Y)}\left(\cdot,\cdot\right).$
    \end{description}
\vspace{0.3cm} 
   
 Assumptions ${\bf H.1}$ and ${\bf H.3}$ are  classical in a nonparametric estimation field and they are the same as those classically used in the nonrecursive case.  The first is satisfied by Gaussian and Eipanechnikov kernels,  while the latter is  checked by linear  processes,   as soon as $f$ is bounded. Note that   ${\bf H.1}$(i)-(ii) are technical conditions, the first  allows the cancellation of the first-order term of Taylor development in the computation of the bias term, while the latter ensures the existence of  the second-order term. 
Much more should be said about assumption  ${\bf H.2. }$ This latter  is particular to the recursive problem and is  clearly  unrestrictive since the choice $h_n = C_nn^{-\nu}$, with $C_n\downarrow c > 0$,  and $0<\nu<1$ is a typical example of bandwidth satisfying {\bf H.2.}  Concerning H.4, the condition {\bf H.4(ii)} is clearly checked if $m$ is a bounded function, and implies that $$\text{E}\left( \max\limits_{1\leq i\leq n} |m(Y_i)|^p \right)= O\left((\ln n)^{p/\theta}\right), \text{ for all } p \geq 1,n \geq 2. $$

The  earlier  condition  was used by Bosq and Cheze-Payaud \cite{bosq-cheze} to study the mean square error of the Nadaraya-Watson  estimator. Assumption {\bf H.4(iii)} was used by Roussas and Tran \cite{Roussas-Tran}  to study the asymptotic normality of $r_n^{DW}({x})$.  
Now, we can give the main result.
  \subsection{Main result}
 Let us set  $$B_n=h_n^2\frac{\beta_{d(1-\ell)+2}}{\beta_{d(1-\ell)}}\frac{1}{2} \sum\limits_{1\leq i,j\leq d}\left(\frac{\partial^2r(x)}{\partial x_i\partial x_j} +2\frac{\partial\ln f(x)}{\partial x_i} \frac{\partial r(x)}{\partial x_j}\right)\int_{\mathbb{R}^d}v_iv_jK(v)dv.$$

The pointwise asymptotic gaussian distribution for our class of nonparametric recursive regression estimate is given in Theorem \ref{normalite_regression1}  below, and will be proved in Section \ref{preuves}.\\
\begin{thrm}\label{normalite_regression1}
When the assumptions    $\textbf{H.1}-\textbf{H.4},$ hold,   if for  all $p>0$, $\left(\ln n\right)^\frac{1}{\theta}h_n^p\rightarrow0, $ as  $n\rightarrow\infty, $   then 
 $$\sqrt{nh_n^d}\left[r_n^\ell(x)-r(x)-B_n\right]\stackrel{\mathcal{L}}\rightarrow \mathcal{N}\left[0,~~ \frac{\sigma_\ell^2(x)V(x)}{f^2(x)}\right], \text{ as  }n\rightarrow\infty,$$
  for all $x$ such that $f(x)>0$,  where
   $$
   \sigma_\ell^2(x)=\frac{\beta_{d(1-2\ell)}}{\beta_{d(1-\ell)}^2} f(x)\int_{\R^d}K^2(x)dx\text{ and } V(x)=\text{E}\left[m^2(Y_0)|X_0=x\right]-r^2(x).\\
$$ 
 \end{thrm}
 One may derive a simpler version of Theorem \ref{normalite_regression1} by  using an additional assumption which allows the cancellation of the bias term $B_n$.
\begin{coro}\label{corollaire} Under  the assumptions    $\textbf{H.1}-\textbf{H.4}$ and  if $nh_n^{d+4}\rightarrow 0$ as $n\rightarrow 0,$ then  
$$\sqrt{nh_n^d}\left[r_n^\ell(x)-r(x)\right]\stackrel{\mathcal{L}}\rightarrow \mathcal{N}\left[0,~~ \frac{\sigma_\ell^2(x)V(x)}{f^2(x)}\right], \text{ as  }n\rightarrow\infty,$$
  for all $x$ such that $f(x)>0$.
\end{coro}
Corollary \ref{corollaire} is an extension of the Rousssas-Tran's \cite{Roussas-Tran} result on Devroye-Wagner estimate to the general family of recursive estimators $r_n^\ell(x)$ for which the Devroye-Wagner estimate is especial case. The condition $nh_n^{d+4}\rightarrow 0$ as $n\rightarrow 0,$  implies that  $\left(\ln n\right)^\frac{1}{\theta}h_n^p\rightarrow0, $ for all $p>0$, and it is satisfied  by the choice $h_n = C_nn^{-\nu}$, with $C_n\downarrow c > 0$ and 1$/(d+4)<\nu<1/(d+2)$.
Let us  mention that  {\bf H.2(iii)} will play a key role in our methodology, in particular when we prove  the negligibility of some covariance terms for   $0\leq\ell\leq(d-2)^+/2 $, but it is  not necessary if $\ell>1/2$.  Also if $\ell>1/2$,  our results can be established for $\rho>2$.  So, we observe that  the estimators built with `small' values of $\ell$  allow some restrictions  on the smooth parameter  $h_n$ and the strong mixing coefficient. However, as shown in Amiri \cite{Amiri2009}, these estimators are preferable than those built with `large' $\ell$ in terms of small variance criterion. \\ 
In practice, the constants of variance appearing in Theorem \ref{normalite_regression1}   need to be estimated. To this end, one may consider the simple Gaussian  kernel and replace $f(x)$ by $f_n^\ell(x)$. There are many possibilities for constructing a consistent conditional variance estimate, one  may use the functional kernel regression technique.\\
In order to prove Theorem \ref{normalite_regression1}, let us consider using  the following decomposition.
$$r_n^\ell(x)-r(x)=\left[\tilde{r}_n^\ell(x)-r(x)\right]+\left[r_n^\ell(x)-\tilde{r}_n^\ell(x)\right],$$
where $\tilde{r}_n^\ell(x)=\tilde{\varphi}_n^{\ell}(x)/f_n^\ell(x)$,  $\tilde{\varphi}_n^{\ell}(x)$ being a  truncated version of  $\varphi_n^{\ell}(x)$ defined by 
 $$\tilde{\varphi}_n^{\ell}(x)=\frac{1}{\sum\limits_{i=1}^nh_i^{d(1-\ell)}}\sum\limits_{i=1}^n\dfrac{Y_i}{h_i^{d\ell}}\textbf{1}_{\left\{\left|Y_i\right|\leq b_n\right\}}K\left( \dfrac{x- X_i}{h_i}\right),$$
with $b_n$, a sequence of real numbers  which goes  to $+\infty$ as $n\rightarrow\infty.$ Next, we need the following preliminary lemmas.
{\lmm \label{Biais-var-errcadr-iid} When  the assumptions  $\textbf{H.1} \text{ and  }\textbf{H.2} $ hold, then for all $\ell\in[0,1]$ 
 \begin{description}
\item[(a)]  $$h_n^{-4}\left[\text{E}f_n^\ell(x)-f(x)\right]^2\longrightarrow
 \left[\frac{\beta_{d(1-\ell)+2}}{\beta_{d(1-\ell)}}\right]^2 b_f^2(x)\text{ as  }n\rightarrow\infty;$$
  \item[(b)] $$h_n^{-4}\left[\text{E}\varphi_n^\ell(x)-\varphi(x)\right]^2\longrightarrow
 \left[\frac{\beta_{d(1-\ell)+2}}{\beta_{d(1-\ell)}}\right]^2 b_\varphi^2(x)\text{ as  }n\rightarrow\infty,$$
 where, if   $h\in C_d^2(b)$,  we  set  $$\label{b2}b_h(x):=\frac{1}{2} \sum\limits_{1\leq i,j\leq d}\frac{\partial^2h}{\partial x_i\partial x_j} (x)\int_{\mathbb{R}^d}v_iv_jK(v)dv.$$

\item[(c)] Moreover if {\bf H.3} holds, then $$nh_n^d\text{Var}f_n^\ell(x)\longrightarrow\sigma_\ell^2(x), \text{ as  }n\rightarrow\infty,$$
for all $x$ such that $f(x)>0.$

\end{description}}
{\it Proof.} The results (a) and (c) of Lemma \ref{Biais-var-errcadr-iid} are obtained    in Amiri \cite{Amiri2009}, while (b) can be established in the same manner as (a) by substituting  $f$ by $\varphi$. \hfill $\square$

{\lmm\label{LEMM2} When  the assumptions  $\textbf{H.1} -\textbf{H.4} $ hold, then for all $\ell\in[0,1]$ 
 \begin{description}
\item[(a)] $$nh_n^d\text{Var}\tilde{\varphi}_n^\ell(x)\longrightarrow\sigma_\ell^2(x)\left[r^2(x)+V(x)\right],\text{ as  }n\rightarrow\infty.$$
\item[(b)]  $$nh_n^d\text{Cov}\left[f_n^\ell(x),\tilde{\varphi}_n^\ell(x)\right]\rightarrow\sigma_\ell^2(x)r(x) \text{ as } n\rightarrow\infty .$$
\end{description}

}

{\it Proof.}
\begin{description}

\item[(a)] Let us set
 $$V_n^*=\sum_{k=1}^n\text{E}Z_{k,n}^{*2}\text{ where }
 Z_{i,n}^*=W_{n,i}-\text{E}W_{n,i},\text{ with } W_{n,i}:=\dfrac{K\left( \dfrac{x-X_i}{h_i}\right)m\left(Y_i\right)\textbf{1}_{\left\{|m(Y_i)|\leq b_n\right\}}}{h_i^{d\ell}}.$$
The variance of $\tilde{\varphi}_n^\ell(x)$ can be decomposed in variance and covariance terms as
$$\text{Var}\tilde{\varphi}_n^\ell(x)=\frac{1}{n^2h_n^{2d(1-\ell)}B_{n,d(1-\ell)}^2}\left[V_n^*+{ \sum_{k=1}^{n}}_{k\neq k'}\sum_{k'=1}^{n}\text{Cov}\left(Z_{k,n},Z_{k',n} \right)\right].$$
Concerning the variance term one may  write
$$\begin{array}{rl}\frac{V_n^*}{nh_n^{d(1-2\ell)}B_{n,d(1-2\ell)}}=&\frac{nh_n^d}{\left(\sum\limits_{i=1}^nh_i^{d(1-\ell)}\right)^2}\sum\limits_{k=1}^{n}\left\{h_k^{-2d\ell} \text{E}K^2\left(\frac{x-X_0}{h_k} \right) m^2(Y_0)\right.\\
-&\left.h_k^{-2d\ell}\text{E}K^2\left(\frac{x-X_0}{h_k} \right)m^2(Y_0)\textbf{1}_{\left\lbrace\mid m(Y_0)\mid>b_n \right\rbrace}\right.\\
-&\left.\text{E}^2K\left(\frac{x-{X}_0}{h_k} \right)m(Y_0){\bf 1}_{\left\{|m(Y_i)|\leq b_n\right\}}\right\}
=:D_1+D_2+D_3.\\
\end{array}$$
Assumptions   $\textbf{H.4} (ii), (iii)$, the dominated convergence theorem and Bochner's lemma imply that  
$$\int_{\R^d}\frac{1}{h_k^d}K^2\left( \frac{x-u}{h_k}\right) \left[V(u)+r^2(u)\right]f(u)du\rightarrow f(x)\left[V(x)+r^2(x)\right]\left\|K\right\|_2^2,\text{as  }k\rightarrow\infty.$$ 
On account of the above,  assumption  $\textbf{H.2}(ii)$ and the Toeplitz lemma allow to deduce that 
$$D_1=\frac{nh_n^d\sum\limits_{k=1}^{n}\left[h_k^{d(1-2\ell)}\int_{\R^d}\frac{1}{h_k^d}K^2\left( \frac{x-u}{h_k}\right) \left[V(u)+r^2(u)\right]f(u)du\right]}{\left(\sum\limits_{i=1}^nh_i^{d(1-\ell)}\right)^2}\rightarrow \sigma_\ell^2(x) \left[V(x)+r^2(x)\right],$$
$\text{as  }n\rightarrow\infty.$
Concerning the term $D_2,$ if $b_n=\left(\delta\ln n\right)^\frac{1}{\theta}\text{ with } \delta>\frac{2}{\lambda},$ then using the  assumptions {\bf H.2(ii)} and  {\bf H.4(ii)}, with the help of  Markov's inequality, we have 
$$\begin{array}{rl}\left|D_2\right|\leq&\frac{\|K\|_\infty^2\left\{\text{E}m^4(Y_0)P\left(|m(Y_0)|>b_n\right)\right\}^\frac{1}{2}
nh_n^d\sum\limits_{k=1}^{n}h_k^{-2d\ell}}{\left(\sum\limits_{i=1}^nh_i^{d(1-\ell)}\right)^2}\leq \frac{\|K\|_\infty^2\left\{\text{E}m^4(Y_0)P\left(|m(Y_0)|>b_n\right)\right\}^\frac{1}{2}
B_{n,-2d\ell}}{h_n^dB_{n,d(1-\ell)}^2}\\
=&O\left[\frac{\exp\left(-\frac{\lambda b_n^\theta}{2}\right)\left(\ln n\right)^\frac{2}{\theta}B_{n,-2d\ell}}{h_n^{d}B_{n,d(1-\ell)}^2}\right]\rightarrow0,\text{ as }n\rightarrow\infty.\end{array}$$
Next for the last term $D_3$, from  {\bf H.2(i)-(iii)} and the logarithmic choice of $b_n$, one may  write 
\begin{eqnarray*}
|D_3|&\leq& \frac{b_n^2nh_n^{d}}{n^2h_n^{2d(1-\ell)}B_{n,d(1-\ell)}^2}\sum\limits_{k=1}^{n}h_k^{-2d\ell}\left(\text{E}K\left(\frac{x-X_i}{h_i}\right)\right)^2=O\left(h_nb_n^2 \right)\rightarrow 0, \text{ as  }n\rightarrow+\infty .
\end{eqnarray*}
Therefore 
$$V_n^*\sim nh_n^{d(1-2\ell)}\beta_{d(1-2\ell)}f(x) \left[V(x)+r^2(x)\right]\int_{\R^d}K^2(u)du,\text{ as }n\rightarrow\infty.$$
It follows that
$$\frac{nh_n^dV_n^*}{n^2h_n^{2d(1-\ell)}B_{n,d(1-\ell)}^2}\rightarrow\frac{\beta_{d(1-2\ell)}f(x) \left[V(x)+r^2(x)\right]}{\beta_{d(1-\ell)}^2}\int_{\R^d}K^2(u)du,\text{ as }n\rightarrow\infty.$$
Now, let us show that the covariance term of Var$\tilde{\varphi}_n^\ell(x)$ is negligible.  To this end, define a sequence $c_n$ of real numbers tending to zero as $n$ goes to infinity, and  write
$$\begin{array}{rl}
\frac{{ \sum\limits_{k=1}^{n}}_{k\neq k'} \sum\limits_{k'=1}^{n}\text{Cov}\left(Z_{k,n},Z_{k',n} \right)}{n^2h_n^{2d(1-\ell)}B_{n,d(1-\ell)}^2}
\leq&\frac{2\left({\sum\limits_{i=1}^n}_{i>j}\sum\limits_{j=1}^n \left|A_{i,j}\right|1_{\{1\leq i-j\leq c_n\}}
+{\sum\limits_{i=1}^n}_{i>j}\sum\limits_{j=1}^n\left|A_{i,j}\right|1_{\{c_n+1\leq i-j\leq n-1\}}\right)}{\left(\sum\limits_{i=1}^nh_i^{d(1-\ell)}\right)^{2}}\\
\leq&\frac{2\left(\sum\limits_{i=1}^{c_n}\sum\limits_{p=1}^nA_{i+p,p}+\sum\limits_{i=c_{n+1}}^{n-1}\sum\limits_{p=1}^nA_{i+p,p}\right)}{\left(\sum\limits_{i=1}^{n}h_i^{d(1-\ell)}\right)^{2}}:=L_1+L_2,\end{array}$$
where $$A_{i+p,p}=\frac{\left|\text{Cov}\left[K\left(\dfrac{x-X_{i+p}}{h_{i+p}}\right)m\left(Y_{i+p}\right)
\textbf{1}_{\left\{|m\left(Y_{i+p}\right)|\leq b_n\right\}},
K\left(\dfrac{x-X_p}{h_p}\right)m\left(Y_p\right)\textbf{1}_{\left\{|m\left(Y_p\right)|\leq b_n\right\}}\right]\right|}{h_{i+p}^{d\ell}h_p^{d\ell}}.$$
On  one hand, the Billingsley  inequality (see e.g.,  Bosq Blanke \cite{bosq-blanke}) implies that
$$A_{i+p,p}\leq 4b_n^2\alpha_X(k)\left\|K\right\|_\infty^2h_{i+p}^{-d\ell}h_p^{-d\ell},$$
and then, it follows from assumptions {\bf H2(ii)} and  {\bf H.4(iv)}  that
$$\begin{array}{rl}
L_2\leq&\frac{8b_n^2\sum\limits_{k=c_n+1}^{n-1}\sum\limits_{p=1}^{n}\alpha_X(k)
h_{p+k}^{-d\ell}h_p^{-d\ell}}{\left\|K\right\|_\infty^2\left(\sum\limits_{i=1}^nh_i^{d(1-\ell)}\right)^{2}}
\leq  \frac{8b_n^2\gamma \left\|K\right\|_\infty^2\sum\limits_{k=c_n}^{n-1}\sum\limits_{p=1}^{n}k^{-\rho}
h_{p+k}^{-d\ell}h_p^{-d\ell}}{\left(\sum\limits_{i=1}^nh_i^{d(1-\ell)}\right)^{2}}
\leq \frac{ 8b_n^2\gamma \left\|K\right\|_\infty^2\dfrac{h_n^{-2d\ell}c_n^{-\rho+1}}{\rho-1}\sum\limits_{p=1}^{n}\left(\dfrac{h_p}{h_n}\right)^{-d\ell}}{\left(\sum\limits_{i=1}^nh_i^{d(1-\ell)}\right)^{2}}\\
    \leq& \dfrac{8b_n^2\gamma \left\|K\right\|_\infty^2c_n^{1-\rho}B_{n,-d\ell}}{nh_n^{2d}B_{n,d(1-\ell)}^2(\rho-1)}.\end{array}$$
    Hence 
$$nh_n^dL_2=O\left(b_n^2c_n^{1-\rho}h_n^{-d}\right).$$
On  the other hand, regarding   about  $L_1$,  one has 
$$\begin{array}{rl}
A_{i+p,p}=&\left|\int_{\mathbb{R}^d}\int_{\mathbb{R}^d}\int_{\mathbb{R}^{d'}}\int_{\mathbb{R}^{d'}} K\left(\frac{x-s}{h_{i+p}}\right)K\left(\frac{x-t}{h_{p}}\right)\frac{m(u)\textbf{1}_{\left\{\mid m(u)\mid \leq b_n\right\}}m(v)\textbf{1}_{\left\{\mid m(v)\mid \leq b_n\right\}}G_{i+p,p}\left(s,u,t,v\right)}{\left(h_{i+p}h_{p}\right)^{d\ell}
\left(\sum_{i=1}^nh_i^{d(1-\ell)}\right)^2}dsdtdudv\right|\\
\leq &\frac{b_n^2\left(h_{k+p}h_p\right)^{d(1-\ell)}\sup\limits_{|k-k'|\geq1}\sup\limits_{(s,t)\in\R^{2d}}\int_{\mathbb{R}^{d'}}\int_{\mathbb{R}^{d'}}\left| G_{k,k'}\left(s,u,t,v\right)\right|dudv}{\left(\sum\limits_{i=1}^nh_i^{d(1-\ell)}\right)^{2}}
.\end{array}$$
Then  \begin{equation}\label{majorJ1}\begin{array}{rl}
  L_1\leq&\frac{2b_n^2\sum\limits_{k=1}^{c_n}\sum\limits_{p=1}^{n-k}h_{p+k}^{d(1-\ell)}h_p^{d(1-\ell)}\sup\limits_{|k-k'|\geq1}\sup\limits_{(s,t)\in\R^{2d}}\int_{\mathbb{R}^{d'}}\int_{\mathbb{R}^{d'}}\left| G_{k,k'}\left(s,u,t,v\right)\right|dudv}{\left(\sum\limits_{i=1}^nh_i^{d(1-\ell)}\right)^{2}}\\
   \leq&\frac{2b_n^2 c_n \sum\limits_{p=1}^nh_p^{2d(1-\ell)}\sup\limits_{|k-k'|\geq1}\sup\limits_{(s,t)\in\R^{2d}}\int_{\mathbb{R}^{d'}}\int_{\mathbb{R}^{d'}}\left| G_{k,k'}\left(s,u,t,v\right)\right|dudv}{\left(\sum\limits_{i=1}^nh_i^{d(1-\ell)}\right)^{2}}.
 \end{array}\end{equation}
 At this point,  we distinguish two cases according to small and large values of $\ell$.\\
 If $\ell\in\left[\left(\frac{d-2}{2d}\right)^+,1\right],$ then $2d(1-\ell)\leq d+2$ replies  $B_{n, 2d(1-\ell)}\rightarrow\beta_{2d(1-\ell)}<\infty,$
$\text{ as }n\rightarrow\infty,$ because of  $\textbf{H.2}(ii).$ It follows that 

 $$L_1\leq\frac{ 2b_n^2c_n B_{n, 2d(1-\ell)}\sup\limits_{|k-k'|\geq1}\sup\limits_{(s,t)\in\R^{2d}}\int_{\mathbb{R}^{d'}}\int_{\mathbb{R}^{d'}}\left| G_{k,k'}\left(s,u,t,v\right)\right|dudv}{nB_{n,d(1-\ell)}^2},$$
which implies that
$$nh_n^dL_1=O\left(b_n^2c_nh_n^{d}\right).$$
Thus, when   $c_n:= \left\lfloor h_n^{-\frac{2d}{\rho}}\right\rfloor,$ and $b_n=\left(\delta\ln n\right)^\frac{1}{\theta}\text{ with } \delta>\frac{2}{\lambda},$ then 
$$\frac{nh_n^d}{n^2h_n^{2d(1-\ell)}B_{n,d(1-\ell)}^2}{ \sum\limits_{k=1}^{n}}_{k\neq k'} \sum\limits_{k'=1}^{n}\text{Cov}\left(Z_{k,n},Z_{k',n} \right)=O\left(b_n^2h_n^{-\frac{d(2-\rho)}{\rho}}\right)\rightarrow0,\text{ as }n\rightarrow\infty,$$ 
$\text{since }  \rho>2.$ Now, if $d\geq3,~~\ell\in\left[0,\frac{d-2}{2d}\right[,$ then the  term $L_1$  cannot be studied as previously,  because assumption {\bf H.2(ii)} is not satisfied   since $2d(1-\ell)>d+2$.  In this case, let us consider relation    (\ref{majorJ1}) and  choose a real number $\xi$ such that $\dfrac{1}{\rho-1}<\xi\leq\dfrac{2}{d}.$ Let us mention  that $\xi$ exists only if $\rho>\dfrac{d+2}{2}.$
Thus, we have the relation $d(\xi+1)\leq d+2,$  which implies that $B_{n,d(\xi+1)}\rightarrow\beta_{d(\xi+1)}<\infty,\text{ as  }n\rightarrow\infty,\text{ by vertue of  $\textbf{H.2}(ii)$.}$ Next, since $h_n$ decreases on has 
$\sum\limits_{i=1}^nh_i^{d(1-\ell)}\geq h_1^{-d\ell}\sum\limits_{i=1}^nh_i^d.$ It follows that 
 $$\frac{c_nb_n^2\sum\limits_{p=1}^nh_p^{2d(1-\ell)}}{\left(\sum_{i=1}^nh_i^{d(1-\ell)}\right)^{2}}\leq \dfrac{c_nb_n^2h_1^{d(1-\xi-2\ell)}nh_n^{d(\xi+1)}B_{n,d(\xi+1)}}{n^2h_1^{-2d\ell}h_n^{2d}B_{n,d}^2}\leq
 \dfrac{c_nb_n^2h_1^{d(1-\xi)}h_n^{d\xi}B_{n,d(\xi+1)}}{nh_n^dB_{n,d}^2},$$
 because $0\leq \ell<\frac{d-2}{2d}\Rightarrow 1-\xi-2\ell>0, \text{ as long  as  }\xi\leq\frac{2}{d}.~$ Therefore,  from (\ref{majorJ1}) we have 
 $$nh_n^dL_1=O\left(c_nb_n^2h_n^{d\xi}\right).$$
 The choices $c_n:= \left\lfloor h_n^{-\frac{d(\xi+1)}{\rho}}\right\rfloor$ and $b_n=\left(\delta\ln n\right)^\frac{1}{\theta}\text{ with } \delta>\frac{2}{\lambda},$ imply the negligibility of the covariance term. 

\item[(b)]  Let us consider  the  decomposition $$\text{Cov}\left[f_n^\ell(x),\tilde{\varphi}_n^\ell(x)\right]=\left[\sum_{i=1}^nh_i^{d(1-\ell)}\right]^{-2}\left[\sum_{i=1}^n
A_{ii}+{\sum_{i=1}^n}_{i\neq j}\sum_{j=1}^nA_{ij}\right]
:=F_1+F_2.$$
where,  for all integers  $s,t$
$$A_{s,t}:=\text{Cov}\left[ \dfrac{1}{h_s^{d\ell}}K\left( \dfrac{x-X_s}{h_s}\right),\dfrac{m(Y_t)}{h_t^{d\ell}}\textbf{1}_{\left\{\left|m(Y_i)\leq b_n\right|\right\}}K\left(\dfrac{x-X_t}{h_t}\right)\right].$$
Next, we proceed as in the proof of (a) and find
$$nh_n^dF_1\rightarrow\sigma_\ell^2(x)r(x), \text{  and   } nh_n^dF_2\rightarrow0, \text{ as }n\rightarrow\infty.$$
\end{description}
\hfill $\square$

\section{Proof of main result}\label{preuves}

{\it Proof.}  To prove the main result, we  show that  the asymptotic distribution  of the principal  term $\left[\tilde{r}_n^\ell(x)-r(x)\right]$ is normal, while the  residual term $\left[r_n^\ell(x)-\tilde{r}_n^\ell(x)\right]$  is negligible. First, observe that if $b_n=\left(\delta\ln n\right)^\frac{1}{\theta}\text{ with } \delta>\frac{2}{\lambda},$ then for all $\varepsilon>0,$ we have 
$$\begin{array}{rl}P\left(\left|\varphi_n^\ell(x)-\tilde{\varphi}_n^\ell(x)\right|>\varepsilon/\sqrt{nh_n^d}\right)
 \leq P\left(\bigcup\limits_{i=1}^n\left\{|Y_i|>b_n\right\}\right)
 \leq nP\left(|Y_0|>b_n\right)
 \leq\text{E}e^{\lambda|m(Y_0)|^\theta}n^{1-\lambda\delta}.\end{array}$$
So,  for all  $\varepsilon>0$, $\sum\limits_{n=1}^\infty P\left(\left|\varphi_n^\ell(x)-\tilde{\varphi}_n^\ell(x)\right|>\varepsilon/\sqrt{nh_n^d}\right)<\infty,$ and the  Borel-Cantelli lemma implies that 
$$\sqrt{nh_n^d}\left[r_n^\ell(x)-\tilde{r}_n^\ell(x)\right]\rightarrow0 \text{ a.s, as }  n\rightarrow\infty.$$
One may  prove  in the same manner  that $f_n^\ell(x)\rightarrow f(x)$ a.s as $n\rightarrow\infty$. 
Next, we need to show that
 $$\sqrt{nh_n^d}\left[\tilde{r}_n^\ell(x)-B_n-r(x)\right]\stackrel{\mathcal{L}}\rightarrow\mathcal{N}\left[0,~~
 \frac{\beta_{d(1-2\ell)}\|K\|_2^2V(x)}{\beta_{d(1-\ell)}^2f(x)}\right],$$
 as  $n\rightarrow\infty.$  To this end, we  use the following representation
      $$
      \tilde{r}_n^\ell(x)-r(x)-B_n
=\frac{1}{f_n^\ell(x)\text{E}f_n^\ell(x)}
 \left[\begin{array}{c}\text{E}f_n^\ell(x)\\-\text{E}\tilde{\varphi}_n^\ell(x)\end{array}\right]^T \left[\begin{array}{c}\tilde{\varphi}_n^\ell(x)-\text{E}\tilde{\varphi}_n^\ell(x)\\ f_n^\ell(x)-\text{E}f_n^\ell(x)\end{array}\right]+o\left(\frac{1}{\sqrt{nh_n^d}}\right).$$
 Now, applying  the Cramer-Wold  device and remembering that $f_n^\ell(x)\stackrel{a.s}\rightarrow f(x), $ and  $\text{E}f_n^\ell(x)\rightarrow f(x), $ as $n\rightarrow\infty,$ the proof of the Theorem  \ref{normalite_regression1} is straightforward from the following claim:
  $$ \sqrt{nh_n^d}\left[\begin{array}{c}\tilde{\varphi}_n^\ell(x)-\text{E}\tilde{\varphi}_n^\ell(x)\\ f_n^\ell(x)-\text{E}f_n^\ell(x)\end{array}\right]\stackrel{\mathcal{L}}\rightarrow\mathcal{N}_2\left\{0,~~ \sigma_\ell^2(x)\left[\begin{array}{lr}V(x)+r^2(x)&r(x)\\ r(x) &1\end{array}\right]\right\},\text{ as }n\rightarrow\infty.$$
 This last convergence  is equivalent to   
  
\begin{equation} \label{Principal}\sqrt{nh_n^d}\left\{\lambda_1\left[f_n^\ell(x)-\text{E}f_n^\ell(x)\right]
+\lambda_2\left[\tilde{\varphi}_n^\ell(x)-
\text{E}\tilde{\varphi}_n^\ell(x)\right]\right\}\stackrel{\mathcal{L}}\rightarrow \mathcal{N}\left[0, \Sigma_\ell^2(x)\right], \text{ as }n\rightarrow\infty,\end{equation}
for each   $\lambda_1,\lambda_2\in \mathbb{R}$ such that $\lambda_1+\lambda_2\neq0$, where $$\Sigma_\ell^2(x):=\sigma_\ell^2(x)\left\{\lambda_1^2+2\lambda_1\lambda_2r(x)+\lambda_2^2\left[V(x)+r^2(x)\right]\right\}.$$
 Hence, the main result will be completely   proven if  (\ref{Principal})  were  established. To this end, let us set
$$\tilde{\Psi}_{nj}:=\lambda_1\Psi_{nj}
 +\lambda_2\Psi_{nj}',$$
where  $\Psi_{nj}:=\left[\frac{h_n^{d(2\ell-1)}}{n}\right]^\frac{1}{2}\frac{h_j^{-d\ell}}{B_{n,d(1-\ell)}}\left(V_{nj}-\text{E}V_{nj}\right)\text{ and }\Psi_{nj}':=\left[\frac{h_n^{d(2\ell-1)}}{n}\right]^\frac{1}{2}\frac{h_j^{-d\ell}}{B_{n,d(1-\ell)}}\left(W_{nj}-\text{E}W_{nj}\right)
 $
 with   $$V_{nj}:=K\left(\frac{x-X_j}{h_j}\right)\text{ and }W_{nj}:=K\left(\frac{x-X_j}{h_j}\right)m\left(Y_j\right)
 \textbf{1}_{\left\{\left|m(Y_j)\right|\leq b_n\right\}}.$$
 Next, consider  the sequences  $\varsigma_n,\tau_n,$ and  $r_n$ defined as
 $$\tau_n:=\left\lfloor\tau_0\log n\right\rfloor, ~~\varsigma_n:=\left\lfloor\frac{\tau_0\sqrt{nh_n^d}}{(\log n)^{\varsigma_0}}\right\rfloor\text{ and  } r_n:=\left\lfloor\frac{n}{\varsigma_n+\tau_n}\right\rfloor,\text{ with     $\tau_0, \varsigma_0>0.$ }$$
To establish (\ref{Principal}), we use  the  classical  Doob \cite{doob}  methodology,  which consists  of splitting  the  term

 $$\sqrt{nh_n^d}\left\{\lambda_1\left[f_n^\ell(x)-f(x)\right]
+\lambda_2\left[\tilde{\varphi}_n^\ell(x)-\varphi(x)\right]\right\}$$ into large blocks separated by  small blocks defined by 

 $$\begin{array}{ll}
     T_{nm}=\sum\limits_{j=k_m}^{k_m+\varsigma_n-1}\tilde{\Psi}_{nj}\text{ (large blocks) }, & T_{nm}'=\sum\limits_{j=l_m}^{l_m+\tau_n-1}\tilde{\Psi}_{nj} \text{ (small blocks) },\\
     \phantom{1}\\
      T_{nr_n+1}'=\sum\limits_{j=\bar{N}+1}^n\tilde{\Psi}_{nj} \text{ (rest of term)}, &   \end{array}
 $$
where $\bar{N}:=r_n(\tau_n+\varsigma_n)$,
and  $\text{ for  } m=1,\ldots,r_n, k_m:=(m-1)(\varsigma_n+\tau_n)+1,~~l_m:=(m-1)(\varsigma_n+\tau_n)+\varsigma_n+1.$
Next, let  us define the partial sums 
$$S_{n1}=\sum\limits_{m=1}^{r_n}T_{nm}, ~~S_{n2}=\sum\limits_{m=1}^{r_n}T'_{nm} \text{ and  }  S_{n3}=T_{nr_n+1}'.$$
Thus, we can write 
$$\sqrt{nh_n^d}\left\{\lambda_1\left[f_n^\ell(x)-f(x)\right]
+\lambda_2\left[\tilde{\varphi}_n^\ell(x)-\varphi(x)\right]\right\}=S_{n1}+S_{n2}+S_{n3} .$$
The goal is to prove that, $\text{E}S_{n2}^2$ and  $\text{E}S_{n3}^2$  converge to zero,   while the asymptotic distribution of $S_{n1}$ is normal. First, observe that 
  \begin{equation}\label{Sn22}\begin{array}{rl}
 \text{E}S_{n2}^2
 =&\sum\limits_{m=1}^{r_n}\text{Var}(T_{nm}')+2\sum\limits_{1\leq i<j\leq r_n}\text{Cov}(T_{ni}',T_{nj}')\\
 =&\sum\limits_{m=1}^{r_n}\sum\limits_{i=l_m}^{l_m+\tau_n-1}\text{Var}\tilde{\Psi}_{ni}+2\sum\limits_{m=1}^{r_n}\sum\limits_{l_m\leq i< j\leq l_m+\tau_n-1}\text{Cov}\left(\tilde{\Psi}_{ni},\tilde{\Psi}_{nj}\right) \\
 +&2\sum\limits_{1\leq i<j\leq r_n}\sum\limits_{s=l_i}^{l_i+\tau_n-1}\sum\limits_{t=l_j}^{l_j+\tau_n-1}\text{Cov}\left(\tilde{\Psi}_{ns},\tilde{\Psi}_{nt}\right)
  :=\Delta_1+\Delta_2+\Delta_3.\end{array}\end{equation}
  The first term in  (\ref{Sn22}), is decomposed as
  $$\begin{array}{rl}\Delta_1=&\sum\limits_{m=1}^{r_n}\sum\limits_{i=l_m}^{l_m+\tau_n-1}\left[\lambda_1^2\text{Var}\Psi_{ni}
 +\lambda_2^2\text{Var}\Psi_{ni}'+2\lambda_1\lambda_2\text{Cov}\left(\Psi_{ni},\Psi_{ni}'\right)\right]
 :=\Delta_{11}+\Delta_{12}+\Delta_{13}
 .\end{array}$$
 Since  $h_n$  decreases, the   choice of  $b_n=\left(\delta\ln n\right)^\frac{1}{\theta}\text{ with } \delta>\frac{2}{\lambda},$ and $\theta>1/\varsigma_0$ with the help of ${\bf H2(iii)}$  implies that 
 $$\begin{array}{rl}\Delta_{11}+\Delta_{12}=&\frac{h_n^{d(2\ell-1)}}{nB_{n,d(1-\ell)}^2}\sum\limits_{m=1}^{r_n}\sum\limits_{j=l_m}^{l_m+\tau_n-1}h_j^{-2d\ell}\left[
\lambda_1^2 \text{Var}K\left(\frac{x-X_j}{h_j}\right)+\lambda_2^2
 \text{Var}K\left(\frac{x-X_j}{h_j}\right)Y_j\textbf{1}_{\left\{\left|m(Y_j)\right|\leq b_n\right\}}\right]\\
 \leq&\frac{r_n\tau_n\left(1+b_n^2\right)\left\|K\right\|_\infty^2\max(\lambda_1^2,\lambda_2^2)}{nh_n^dB_{n,d(1-\ell)}^2}\rightarrow 0,\text{ as  }n\rightarrow\infty,\end{array}$$
and, similarly,   we have $\Delta_{13}\leq\frac{2\lambda_1\lambda_2b_nr_n\tau_n\left\|K\right\|_\infty^2}
{nh_n^dB_{n,d(1-\ell)}^2}\rightarrow 0,\text{ as } n\rightarrow\infty.$
In the same manner and  by also using the  Cauchy-Schwartz's inequality, we get $$\Delta_2\leq\frac{r_n\tau_n^2\left(1+b_n\right)^2\left\|K\right\|_\infty^2\max(\lambda_1^2,\lambda_2^2)}{nh_n^dB_{n,d(1-\ell)}^2}
 \rightarrow 0,\text{ as  }n\rightarrow\infty.$$
 The last term in  (\ref{Sn22}) is bounded  by Billingsley inequality with the help of assumptions ${\bf H2(iii) }$ and  $\textbf{H.3}(i)$, as follows.
  $$\begin{array}{rl}\Delta_3=&2\sum\limits_{1\leq i<j\leq r_n}\sum\limits_{s=l_i}^{l_i+\tau_n-1}\sum\limits_{t=l_j}^{l_j+\tau_n-1}\left\{\lambda_1^2\text{Cov}\left(\Psi_{ns},\Psi_{nt}\right)\right.
 +\lambda_2^2\text{Cov}\left(\Psi_{ns}',\Psi_{nt}'\right)\\
 +&\left.\lambda_1\lambda_2\left[\text{Cov}\left(\Psi_{ns},\Psi_{nt}'\right)
 +\text{Cov}\left(\Psi_{nt},\Psi_{ns}'\right)\right]\right\}\\
 \leq&\frac{2\left(1+b_n\right)^2\left\|K\right\|_\infty^2\max(\lambda_1^2,\lambda_2^2)h_n^{d(2\ell-1)}}{nB_{n,d(1-\ell)}^2}\sum\limits_{k=1}^{r_n-1}\sum\limits_{j=1}^{r_n}\sum\limits_{s=l_j}^{l_j+\tau_n-1}
\sum\limits_{t=l_j}^{l_j+\tau_n-1}(h_sh_t)^{-d\ell}\alpha_X\left[k\left(\varsigma_n+\tau_n\right)\right]\\
\leq&\frac{2\gamma\left(1+b_n\right)^2\left\|K\right\|_\infty^2\max(\lambda_1^2,\lambda_2^2)r_n\tau_n^2}{nh_n^dB_{n,d(1-\ell)}^2}\sum\limits_{k=1}^{r_n-1}e^{-\rho k\tau_n}.\end{array}$$
 Therefore,
 $$\Delta_3=O\left\{\frac{b_n^2r_n\tau_n^2e^{-\rho\tau_n} }{nh_n^{d}B_{n,d(1-\ell)}^2}\left[1-e^{-\rho\tau_n(r_n-1)}\right]\right\}\rightarrow 0, \text{ as } n\rightarrow\infty,$$
 as long as $b_n=\left(\delta\ln n\right)^\frac{1}{\theta}\text{ with } \delta>\frac{2}{\lambda},$ and $\theta>1/\varsigma_0$. Now, let us prove that  $\text{E}S_{n3}^2\rightarrow 0$ as  $n\rightarrow0$. One has
\begin{equation}\label{ES3}\text{E}S_{n3}^2=\sum_{j=\bar{N}+1}^n\text{Var}\tilde{\Psi}_{nj}
+2\sum_{\bar{N}
+1\leq i<j\leq n}\text{Cov}(\tilde{\Psi}_{ni},\tilde{\Psi}_{nj})
:=\Theta_{n1}+\Theta_{n2}.\end{equation}
The  variance term  $\Theta_{n1}$ may be written as 
$$
\Theta_{n1}=\sum_{j=\bar{N}+1}^n\left[\lambda_1^2\text{Var}\Psi_{nj}
+\lambda_2^2\text{Var}\Psi_{nj}'+2\lambda_1\lambda_2\text{Cov}\left(\Psi_{nj},\Psi_{nj}'\right)\right]
:=\lambda_1^2\Theta_{n11}+\lambda_2^2\Theta_{n12}+2\lambda_1\lambda_2\Theta_{n13}.
$$
The first term of the  right hand side of the preview decomposition  satisfies  the relation
$$nh_n^d\text{Var}f_n^\ell(x)\sim\sum_{j=1}^n\text{Var}(\Psi_{nj})=\sum_{j=1}^{\bar{N}}\text{Var}
(\Psi_{nj})+\Theta_{n11}.$$
However,  one may write 
  $$\sum_{j=1}^{\bar{N}}\text{Var}
(\Psi_{nj})=\left(\frac{nh_n^d}{\bar{N}h_{\bar{N}}^d}\right)\bar{N}
h_{\bar{N}}^d\text{ Var } f_{\bar{N}}^\ell(x).$$ 
Since  $\bar{N}\sim n,$  the condition $u_n\sim v_n$ implies  $h_{u_n}\sim h_{v_n},$  which leads to $nh_n^d\sim\bar{N}h_{\bar{N}}^d,$ which together with Lemma  \ref{Biais-var-errcadr-iid}(c) imply  that 
 $\sum\limits_{j=1}^{\bar{N}}\text{Var}(\Psi_{nj})\rightarrow\sigma_\ell^2(x),\text{ as }n\rightarrow\infty.$
  It follows that $\Theta_{n11}=o(1),$  because $\sum\limits_{j=1}^n\text{Var}(\Psi_{nj})\rightarrow\sigma_\ell^2(x),\text{ as }n\rightarrow\infty.$
Let us mention that  if  $\ell\geq1/2$, then  the condition $u_n\sim v_n$ implies  $h_{u_n}\sim h_{v_n},$ is not necessary. Indeed, the variance term   $\Theta_{n1}$ can be written as 

$$\Theta_{n1}=\frac{B_{n,d(1-\ell)}^{-2}}{n}\sum_{i=\bar{N}+1}^n
\left(\frac{h_i}{h_n}\right)^{d(1-2\ell)}h_i^{-d}\text{Var}\left\{K\left(\frac{x-X_i}{h_i}\right)\left(1+m(Y_i){\bf 1}_{|m(Y_i)|\leq b_n}\right)\right\}.$$
Since  $h_n$ is decreasing and $\ell\geq\frac{1}{2}$, then  the Toeplitz lemma,  with the help of  assumption {\bf H.2(ii)}  and the convergence $$h_i^{-d}\text{Var}K\left(\frac{x-X_i}{h_i}\right)\rightarrow f(x)\int_{\R^d}K^2(x)dx,\text{ as  }i\rightarrow\infty$$ imply  that 
   $$\Theta_{n1}\leq\frac{\text{Cste}\left(n-\bar{N}\right)(1+b_n)^2}{nB_{n,d(1-\ell)}^{2}}.$$
  Because of  $n-\bar{N}\leq \varsigma_n+\tau_n,$ it follows that  $\Theta_{n1}\rightarrow0$ as $n\rightarrow\infty$, provided $b_n=\left(\delta\ln n\right)^\frac{1}{\theta}\text{ with } \delta>\frac{2}{\lambda},$ and $\theta>1/\varsigma_0$.
Also,  in the same manner,  and  by replacing $f_n^\ell$ by $\tilde{\varphi}_n^\ell$,  we can deduce from Lemma \ref{LEMM2}(a)   that   $\Theta_{n12}=o(1)$. 
Finally, the last term  $\Theta_{n13}$ is bounded  similarly to the first term  by using Lemma \ref{LEMM2}(b).
Therefore, from $\Theta_{n13}=o(1),$ it follows that 
$\Theta_{n1}\rightarrow0$ as  $n\rightarrow\infty$.  Now, let us study the term  $\Theta_{n2}$ in  (\ref{ES3}). This can be decomposed as 
$$\Theta_{n2}=2\sum_{\bar{N}
+1\leq i<j\leq n}\left[\lambda_1^2\text{Cov}\left(\Psi_{ni},\Psi_{nj}\right)
+\lambda_2^2\text{Cov}\left(\Psi_{ni}',\Psi_{nj}'\right) +2\lambda_1\lambda_2
\text{Cov}\left(\Psi_{ni},\Psi_{nj}'\right)\right].$$
As in  the proof of  Lemma \ref{LEMM2}, one may  show that
 $$\sum\limits_{\bar{N}
+1\leq i<j\leq n}\left[\lambda_1^2\text{Cov}\left(\Psi_{ni},\Psi_{nj}\right)
+\lambda_2^2\text{Cov}\left(\Psi_{ni}',\Psi_{nj}'\right)\right]\rightarrow0,\text{ as } \rightarrow0$$ and 
$$\sum\limits_{\bar{N}+\leq i<j\leq n}\text{Cov}\left(\Psi_{ni},\Psi_{nj}'\right)\leq\sum\limits_{1\leq i<j\leq n}\text{Cov}\left(\Psi_{ni},\Psi_{nj}'\right)\rightarrow0,\text{ as } \rightarrow0.$$
Hence,  $\Theta_{n2}\rightarrow0,\text{ as } \rightarrow0.$ 
 To complete the  proof we must show that  the asymptotic distribution of $S_{n1}$ is  normal.  To this end let us check the Lindeberg-Feller conditions for $S_{n1}$. First,  we  consider a sequence of  iid random variables $Z_{n1},\ldots,Z_{nr_n}$,   having the same distribution as  $T_{nm}$.  Then,   $\text{E}Z_{n1}=0$ and  if    $\Phi_{T_{nm}}$ is the characteristic function   (f.c.) of   $T_{nm}$, then   $\Phi_{T_{nm}}^{r_n}$ is the f.c.  of the random  variable  $\sum\limits_{m=1}^{r_n}Z_{nm}$.  To establish the asymptotic normality of   $S_{n1}$, it suffices  to prove that the  variables   $\sum\limits_{m=1}^{r_n}Z_{nm}$ and   $\sum\limits_{m=1}^{r_n}T_{nm}$  have the same distribution, and that this latter is Gaussian.  By the Volkonskii-Rosanov \cite{volkonsky-rozanov} lemma, one has 
$$\left|\text{E}\prod_{m=1}^{r_n}e^{itT_{nm}}-\prod_{m=1}^{r_n}\text{E}e^{itT_{nm}}\right|\leq8(r_n-1)\alpha(\tau_n)\leq \rho_0r_ne^{-\rho_1\tau_n}
\rightarrow0,\text{ as } n\rightarrow\infty.$$
It follows that $\left|\text{E}\prod_{m=1}^{r_n}e^{itT_{nm}}-\Phi_{T_n}^{r_n}\right|\rightarrow0, \text{ as  }n\rightarrow\infty.$ Then, it suffices to prove that $\Phi_{T_{nm}}^{r_n}$ converges to the  characteristic function of a Gaussian random variable. To this end, we proceed as follows. Set 
$Z_{nm}':=\frac{Z_{nm}}{s_n},\text{ where  }s_n^2:=\sum_{m=1}^{r_n}\text{Var}Z_{nm}.$
One has $s_n^2\rightarrow\Sigma_\ell^2(x),\text{ as  } n\rightarrow\infty.$
Indeed, $s_n^2=\sum\limits_{m=1}^{r_n}\text{Var}T_{nm}\rightarrow\Sigma_\ell^2(x),\text{ as  } n\rightarrow\infty,$  because on  one hand we have from Lemmas \ref{Biais-var-errcadr-iid} and  \ref{LEMM2}:
$$\text{Var}S_{n1}\sim nh_n^d\left\{\lambda_1^2\text{Var}f_n^\ell(x)+\lambda_2^2\text{Var}\tilde{\varphi}_n^\ell(x)+2\lambda_1\lambda_2\text{Cov}\left[ f_n^\ell(x),\tilde{\varphi}_n^\ell(x)\right]\right\}\rightarrow\Sigma_\ell^2(x),\text{ as }n\rightarrow\infty,$$
 and one may show, on  the other hand, as for   $\Delta_2$,  that $\sum\limits_{1\leq i<j\leq r_n}\text{Cov}(T_{ni},T_{nj})\rightarrow0,\text{ as }n\rightarrow\infty.$
Hence, the   variables $Z_{nm}'$ are  iid, $\text{E}Z_{n1}'=0\text{ and  }\sum\limits_{m=1}^{r_n}\text{Var}Z_{nm}'=1.$ 
By virtue of the Lindeberg conditions (c.f. Lo\`eve \cite{loeve}), we have to show that for all  $\varepsilon>0,$
$$\sum\limits_{m=1}^{r_n}\text{E}\left(Z_{nm}'^2\textbf{1}_{\left\{\left|Z_{nm}'\right|>\varepsilon\right\}}\right)\rightarrow0,\text{ as }n\rightarrow\infty.$$
Noting that 
$\left|T_{nm}\right|\leq
\frac{\varsigma_n\left\|K\right\|_\infty\left(1+b_n\right)}{\sqrt{nh_n^d}B_{n,d(1-\ell)}},$
and applying  Markov's inequality, one has 
 $$\begin{array}{rl}
\sum\limits_{m=1}^{r_n}\text{E}\left(Z_{nm}'^2\textbf{1}_{\left\{\left|Z_{nm}'\right|>\varepsilon\right\}}\right)=&\sum\limits_{m=1}^{r_n}\text{E}
\left(\frac{T_{nm}^2}{s_n^2}
\textbf{1}_{\left\{\left|T_{nm}\right|>\varepsilon s_n\right\}}\right)
\leq\frac{\varsigma_n^2\left(1+b_n\right)^2\left\|K\right\|_\infty^2}{nh_n^{d}B_{n,d(1-\ell)}^2s_n^2}\sum\limits_{m=1}^{r_n}P\left( \left|T_{nm}\right|>\varepsilon s_n\right)\\
\leq&\left[\frac{\varsigma_n\left(1+b_n\right)}{\sqrt{nh_n^d}}\cdot\frac{\left\|K\right\|_\infty\varepsilon^{-1}}{s_nB_{n,d(1-\ell)}}\right]^2\rightarrow0, \text{ as }n\rightarrow\infty, \end{array}$$
if $b_n=\left(\delta\ln n\right)^\frac{1}{\theta}\text{ with } \delta>\frac{2}{\lambda},$ and $\theta>1/\varsigma_0$.
\hfill $\square$






\begin{thebibliography}{00}
 \bibitem{Ahmad-lin}
 Ahmad, I. and  Lin, P.E. (1976).  Nonparametric sequential estimation of a multiple regression function, {\it  Bull. Math. Statist.  17,   63-75.}
 \bibitem{Amiri2009}
 Amiri, A. (2009). Sur une famille param\'etrique d'estimateurs s\'equentiels de la densit\'e pour un processus fortement m\'elangeant, {\it C. R. Acad. Sci. Paris, Ser. I, 347, 309-314}.

\bibitem{Amiri} Amiri, A. (2012). Recursive regression estimators with application to nonparametric prediction. \textit{J. Nonparam. Statist.}, \textbf{24}, 169-186.
 \bibitem{bosq-blanke}
Bosq, D. and   Blanke, D. (2007). Inference and prediction in large dimensions,  {\it Wiley series in probability and statistics.}
\bibitem{bosq-cheze}
Bosq, D. and Cheze-Payaud, N. (1999). Optimal asymptotic quadratic error of nonparametric regression function estimates for a continuous-time process from sampled-data, {\it Statistics 32 (3), 229-247.}
\bibitem{devroye81}
Devroye, L. (1981). On the almost everywhere convergence of nonparametric regression function estimates, {\it Ann.  Statist. 9, 1310-1319.}
 \bibitem{dev-wag}
Devroye, L. and  Wagner, T. J. (1980).  On the $L^1$ convergence of kernel estimators of regression functions with application in discrimination, {\it Z. Wahrschein. Verw. Get. 51, 15-25.}
\bibitem{doob}
Doob, J. (1953). Stochastic process, {\it  Wiley New York.}
\bibitem{Grebilicki-pawlak}
Greblecki, W. and Pawlak, M. (1987). Necessary and sufficient consistency conditions for a recursive kernel regression estimate,  {\it J. Multivariate Anal.   23, 67-76.}
\bibitem{GyorfiHardle}
Gyorfi, L., Hardle, W., Sarda, P. and  Vieu, P. (1989). Nonparametric curve estimation from time Series, {\it  Lecture notes in statistics,  Springer}.
\bibitem{GyorfiKholer}
Gy\"orfi, L.,  Kh\"oler, M., Krzy$\dot{\text{z}}$ak, A.  and Walk, H. (2002). A distribution-free theory of nonparametric regression, {\it Springer-Verlag New york.}
 
\bibitem{KRZ}
Krzy$\dot{\text{z}}$ak, A. (1992). Global convergence of the recursive kernel regression estimates with applications in classification and nonlinear system estimation, {\it IEEE Trans. Inform. Theory  38, 1323-1338.}
\bibitem{loeve}
Lo\`eve, M. (1963).  Probability theory, {\it Princeton, New Jersey Van Nostrand.}
\bibitem{prakassa-rao1}
Prakasa-Rao, B.L.S. (1983), Nonparametric functional estimation, {\it New-York: Academic Press.}
\bibitem{Qin}
Qin, Y.S. (1995), ÔAsymptotic distribution of a recursive kernel estimator for a nonparametric regression function
under dependent samplingÕ, {\it Mathematica Applicata, 8(1), 7--13.}
\bibitem{Roussas90}
Roussas, G.G. (1990). Nonparametric regression estimation under mixing conditions.  {\it Stochastic Process. Appl.  36 (1), 107--116.}
\bibitem{Roussas-Tran}
Roussas, G.G. and Tran, L.T. (1992).  Asymptotic normality of the recursive kernel regression estimate under dependence conditions, {\it Annals of Statist. 20 (1), 98-120.}
\bibitem{volkonsky-rozanov}
Volkonskii, V.A. and   Rozanov, Yu.A. (1959). Some limit theorems for random
functions, {\it Theory Probab. Appl. 4,  178-197.}
\bibitem{walk}
Walk, H. (2001). Strong universal pointwise consistency of recursive regression estimates,  {\it Ann. Inst. Statist. Math.  53 (4),  691-707.}
\bibitem{Wang}
Wang, L., and Liang, H.Y. (2004), ÔStrong uniform convergence of the recursive regression estimator under $\varphi$-Mixing
ConditionsÕ, {\it Metrika, 59(3), 245--261.}
\bibitem{Yoshihara}
Yoshihara, K. (1994), Weakly dependent stochastic sequences and their applications: curve estimation based on weakly dependent data {\it (Vol. IV), Tokyo: Sanseido.}

 \end{thebibliography}







\end{document}